\documentclass[11pt]{amsart}
\usepackage{graphicx}
\usepackage{latexsym}
\usepackage{amsfonts,amsmath,amssymb}
\usepackage{algorithm}
\newtheorem{theorem}{Theorem}[section]

\usepackage{xcolor}

\def\eps{\varepsilon}
\def\la{\lambda}
\def\a{\alpha}
\def\be{\beta}

\def\part{\partial}

\newcommand{\beq}{\begin{equation}}
\newcommand{\eeq}{\end{equation}}

\theoremstyle{remark}

\numberwithin{equation}{section}
\date{\today}
\begin{document}
\title[Perfect matching]{On a perfect matching in a random bipartite digraph with average out-degree below two.}
\date{}
{\author{Michal Karo\'nski, Ed Overman and Boris Pittel}}

\begin{abstract} Existence of a perfect matching in a random bipartite digraph with bipartition $(V_1, V_2)$, $|V_i|=n$, 
is studied. The graph is generated in two rounds of random selections of a potential matching partner such that the average
number of selections made by each vertex overall is below $2$. More precisely,
in {the} first round each vertex chooses a potential mate uniformly
at random, and independently of all vertices. Given a fixed integer $m$, a vertex is classified as unpopular if it has been chosen
by at most $m$ vertices from the other side. Each unpopular vertex makes yet another uniform/independent selection of a potential mate. 
The expected number of selections made by a generic vertex $v$, i.e. its out-degree, is asymptotic to $1+\Bbb P(\text{Poisson}(1)\le m)\in (1,2)$.
Aided by Matlab software, we prove that for $m=1$,  whence for all $m\ge 1$, the resulting bipartite graph has a perfect matching a.a.s.\@ (asymptotically almost surely).
On the other hand, for $m=0$ a.a.s.\@ a perfect matching does not exist. This is a thorough revision of the  joint paper 
(JCT(B) 88 (2003), 1-16) by the first author and the third author.
  
\end{abstract}

\keywords
{bipartite graphs, perfect matchings, random, asymptotics}
\subjclass[2010] {05C05, 05C07, 05C30, 05C80, 60C05}

\maketitle
\section{Introduction and main result} A standard model $B_n(d)$ of a random bipartite (di)graph with bipartition $(V_1, V_2)$, $|V_i|=n$, is
generated by each vertex $v\in V_1\cup V_2$ making $d$ uniformly random, independent selections of a potential match from the other side. By computing the
expected number of perfect matchings, Walkup \cite{Wal} proved that asymptotically almost surely (a.a.s.\@) the graph $B_n(1)$ has no perfect
matching. In fact, Meir and Moon \cite{MeiMoo} (cf. Frieze \cite{Fri}) earlier proved that the maximum matching number of $B_n(1)$ is a.a.s.\@
about $0.866 n$. It is not much more difficult to show, using Hall's Marriage Lemma, that a.a.s.\@ the graph $B_n(3)$ does have a perfect matching.
Remarkably, Walkup managed to show that a.a.s.\@ so does the graph $B_n(2)$. Frieze \cite{Fri} was able to prove an analogous result for a non-bipartite graph using Tutte's criterion. 

In this paper we study existence of a perfect matching in a bipartite random graph $B_{n,m}$ which is sandwiched between
$B_n(1)$ and $B_n(2)$. $B_{n,m}$ is generated in two rounds of random selections of a potential match by every vertex $v\in V_1\cup V_2$.
Specifically, in the first round each vertex selects a vertex from the opposite side uniformly at random, and independently of all other vertices. {\it We call
a vertex ``unpopular'' if it has been selected by at most $m$ vertices.\/}  (This definition depends on the value $m$: the larger $m$ the larger the set
of unpopular vertices.) Each unpopular vertex makes yet another uniformly random, and independent selection of a vertex from the other side. (If $m=\infty$, then effectively all the vertices select uniformly at random and independently two vertices from the other side, so $B_{n,\infty}=B_n(2)$.)  The number of vertices that have selected a given vertex is distributed binomially with $n$ trials and success probability $1/n$;
thus it is $\text{Poisson}(1)$ in the limit. It follows that the expected out-degree of a generic vertex in $B_{n,m}$ is 
\[
1+\Bbb P(\text{Poisson}(1)\le m)=1+e^{-1}\sum_{j=0}^m\frac{1}{j!}\uparrow 2,\quad m\to\infty.
\]
Thus the average out-degree of a vertex in $B_{n,m}$ is strictly between $1$ and $2$. Loosely, we can interpret $B_{n,m}$ as $B_n(d_m)$,
$d_m:=1+e^{-1}\sum_{j\le m}1/j!$. In the joint paper \cite{KarPit} the first author and the third author stated and gave a proof of the claim: a.a.s.\@ $B_{n.0}$ (i.e. $B_n(1+1/e)$) has a perfect matching. Recently Michael Anastos and Alan Frieze \cite{AnaFri} pointed out a simple oversight in the proof. We realized
that the oversight invalidated the claim.  A thorough revision of the method in \cite{KarPit} has allowed us to prove 
{
\begin{theorem}\label{1} Let $m\ge 1$.
\begin{align}
\tag{\bf 1} 
&\qquad\Bbb P\bigl(B_{n,m}\text{ has a perfect matching}\bigr)\ge 1-O(n^{-c_m+o(1)}),\\
\notag
&c_m:=1-\frac{1.5}{1+(m+1)\left(1+e^{-1}\sum_{j\le m}1/j! -\log 2\right)},
\end{align}
where $c_1\approx 0.514,\,c_m\uparrow 1 (m\to\infty).$
\begin{equation}
\tag{\bf 2} 
\Bbb P(B_{n,m}\text{ is connected})=1-{O\bigl(n^{-c_m+o(1)}\bigr)}.
\end{equation}
\end{theorem}
}
Let the reader beware that our rigorous proof techniques produced an explicit function $H_{n,m}(t; \bold r)$, $t\in (0,1/2]$, $\bold r\in \Bbb R^4$,
such that, to complete the proof, $\min_{\,\bold r}H_{n,m}(t; \bold r)$ needed to be proved negative for all $t\in (0,1/2]$. This is because the minimum is
the best achievable upper bound for the scaled logarithm of the {\it expected\/}  number of Hall's subgraphs of a given size that are present if there is no
perfect matching. 
By upper-bounding  $\min_{\,\bold r}H_{n,m}(t; \bold r)$, we checked negativity ``manually'' for small $t$, but deferred to Matlab algorithmic software to handle the remaining $t$'s.

As for $m=0$, Matlab revealed that $\min_{\,\bold r} H_{n,0}(t;\bold r)>0$ for all, but very small $t\in (0,1/2]$. 
Of course, positivity of this minimum even for all $t\in (0,1/2]$ does not imply almost sure non-existence of a perfect matching. However this numerical evidence prodded us to try and prove that, contrary to our long-held belief, existence of a perfect matching is indeed highly unlikely.
Using the necessity part of Hall's Lemma we prove the opposite of the claim in \cite{KarPit}:   
\begin{theorem}\label{2} 
\[
\Bbb P\bigl(B_{n,0}\text{ has no perfect matching}\bigr)= 1-O(n^{-1/2+o(1)}).
\]
\end{theorem}
\noindent For the proof itself,  Matlab assistance was not required.
\section{Proof of Theorem \ref{1}.} 
Let $m>0$. {\bf Part 1.\/} If  the graph $B_{n,m}$ has no perfect matching, than by Hall's Marriage Lemma there exist a set $K$ of row vertices 
(elements of $V_1$), or of
column vertices (elements of $V_2$), such that $|L|<|K|$, where $L=\Gamma(K)$ is the set of neighbors of $K$ in $B_{n,m}$. We call such $(K,L)$ ``bad'' pairs.
We focus on the minimal pairs, minimal in a sense that there is no $K'\subset K$ such that $(K',\Gamma(K'))$ is a bad pair. For a minimal bad 
pair $(K,L)$, it is {\it necessary} that (i) $|L|=|K|-1$; (ii)
{ $|K|\in \bigl[2, \lceil n/2\rceil\bigr]$}; (iii) every vertex in $L$ has at least two neighbors in $K$.
Let
{ $E_{nk}$, $k\in \bigl[2, \lceil n/2 \rceil\bigr]$}, denote the expected number of the minimal bad pairs $(K,L)$, with $|K|=k$. We need to prove that 
 $\sum_k E_{nk}\to 0$. By symmetry,
 \begin{equation}\label{1n}
E_{nk}\le 2\binom{n}{k}\binom{n}{k-1}\mathcal P_{nk};
\end{equation}
here $\mathcal P_{nk}$ is the probability that $K=[k]=\{1,\dots,k\}\subset V_1$ and $L=[k-1]=\{1,\dots,k-1\}\subset V_2$ satisfy the conditions
(ii) and (iii). In fact we {\it broaden\/} the condition (iii) a bit, replacing it with (iii'): in two rounds of selections every vertex in $L$ at least twice either selected  a vertex in $K$ or was selected by a vertex in $K$.

\subsection{Case $k\le n^{1/2}$.} On the event in question, let $Y$ ($X$ resp.) be the number of columns in $[k-1]$ (the number of rows in $[n-k]$ resp.)
that selected rows in $[k]$ (columns in $[n-k+1]$ resp.) in the first round. Then the number of unpopular rows in $[k]$ (unpopular columns in $[n-k+1]$
resp.) is at least $k-Y/(m+1)$ ($(n-k+1)-X/(m+1)$ resp.). { (Indeed, every popular row in $[k]$ is selected by at least $m+1$ columns out of $Y$
columns.)} $Y$ and $X$ are independent binomials with parameters $(k-1, k/n)$ and $(n-k, (n-k+1)/n)$
respectively. Therefore
\begin{multline*}
\mathcal P_{nk}\le \left(\frac{k-1}{n}\right)^k \left(1-\frac{k}{n}\right)^{n-k+1}\\
 \times\sum_{j=0}^{k-1}\binom{k-1}{j}\!\left(\frac{k}{n}\right)^j\!\left(1-\frac{k}{n}\right)^{k-1-j}\! \left(\frac{k-1}{n}\right)^{k-j/(m+1)}\\
\times\sum_{i=0}^{n-k}\binom{n-k}{i} \!\left(\frac{n-k+1}{n}\right)^i\!\left(1-\frac{n-k+1}{n}\right)^{n-k-i}\! \left(1-\frac{k}{n}\right)^{n-k+1-i/(m+1)}.
\end{multline*}
{\it Explanation\/} The first line of the RHS is the probability that the first round choices made by rows from $[k]$ (columns from $[n-k+1]$ resp.) are columns from $[k-1]$ (rows from $[n-k]$ resp.). The second line (third line resp.) is an upper bound for the conditional probability that the second round choices made by unpopular rows from $[k]$ (unpopular columns from $[n-k+1]$ resp.) are still columns from $[k-1]$ (rows from $[n-k]$ resp.). Further, the first sum equals
\begin{align*}
&\left(\frac{k-1}{n}\right)^{k}\cdot \left[1+\frac{k}{n}\cdot\left(\frac{k-1}{n}\right)^{-1/(m+1)}\!\!-\frac{k}{n}\right]^{k-1}\\
&\qquad\qquad=\left(\frac{k-1}{n}\right)^{k}\cdot \exp\left(O\Bigl(k^{\frac{m}{m+1}}\Bigr)\right).
\end{align*}
The second sum equals
\begin{multline*}
\left(1-\frac{k}{n}\right)^{n-k+1}\cdot\left[\left(1-\frac{k-1}{n}\right)\left(1-\frac{k}{n}\right)^{-1/(m+1)}+\frac{k-1}{n}\right]^{n-k}\\
=\left(1-\frac{k}{n}\right)^{n-k+1}\left(1+\frac{k}{(m+1)n} +O\bigl(k^2/n^2)\bigr)\right)^{n-k}\\
=\left(1-\frac{k}{n}\right)^{n-k+1}\exp\Bigl(\frac{k}{m+1}+O(k^2/n)\Bigr).
\end{multline*}
Therefore 
\begin{align*}
\mathcal P_{nk}&\le \left(\frac{k-1}{n}\right)^{2k} \left(1-\frac{k}{n}\right)^{2(n-k+1)}\exp\Bigl(\frac{k}{m+1}+O\bigl(k^{\frac{m}{m+1}}+k^2/n\bigr)\Bigr)\\
&=\left(\frac{k}{n}\right)^{2k}\exp\Bigl(-k\frac{2m+1}{m+1}+O\bigl(k^{\frac{m}{m+1}}+k^2/n\bigr)\Bigr).
\end{align*}
Consequently we have
\begin{multline}\label{2n+}
E_{nk}=O\!\left(\!\frac{k}{n}\binom{n}{k}^2\!\left(\frac{k}{n}\right)^{2k}\!\! \exp\Bigl(-k\frac{2m+1}{m+1} +O\bigl(k^{\frac{m}{m+1}}+k^2/n\bigr)\right)\\
=n^{-1} \exp\Bigl(\frac{k}{m+1} +O\bigl(k^{\frac{m}{m+1}}+k^2/n\bigr)\Bigr).
\end{multline}
In particular, for $\eps>0$, $\sum\limits_{k\le (m+1-\eps)\log n}\!\!\!\!\!\! E_{nk} =O(n^{-\eps/(m+1)})$, so that 
\begin{equation}\label{2n++}
\Bbb P\bigl(\exists\text{ a bad pair }(K,L): k\le (m+1-\eps)\log n\bigr)=O(n^{-\eps/(m+1)}).
\end{equation}

\subsection{Case $k\in [n^{1/2}, (n+1)/2]$.}  We write $\mathcal P_{nk}=\Bbb P(A\cap B\cap C)$. 

$A$: first round choices of the rows from $[k]$ are among the columns in $[k-1]$, and for every ``unpopular''  row in $[k]$ (i. e. a receiver of at most
$m$ first-round proposals), its second round choice is still one  of the columns from $[k-1]$.

$B$: first round choices of the columns from $[n-k+1]:=V_2\setminus [k-1]$ are among the rows in $[n-k]:=V_1\setminus [k]$, and, for every
column  $j\in [n-k+1]$ unpopular among the rows in $[n-k]$ in the first round, $j$'s second round choice---in case it is unpopular among rows in $[k]$ too---would still be a row in $[n-k]$.

$C$:  overall, every column vertex from $[k-1]$ has taken part, as a proposer or a ``proposee'', in at least two contacts with the row vertices from $[k]$. 

Clearly the second  round choices of rows from $[n-k]$ are irrelevant for the events $A,\,B,\,$ and $C$. Let $\mathcal G$ denote the (muti)graph, with labeled edges, 
induced by the two rounds of selections  by the row set $[k]$, and the first round selections  by the column set $[n]$. Let $\mathcal H$
be the graph induced by the second round choices by the column set $[n]$. Then $\mathcal G$
is independent of $\bold X=(X_1,\dots,X_n)$, where $X_j$ is the number of first round selections of column $j$ by rows in $[n-k]$, and 
the distribution of $\mathcal H$ conditioned on $\{\bold X=\bold x,\,\mathcal G=G\}$ is the same no matter what the marginal distribution of $\bold X$
is. In the selection process $\bold X$ is distributed multinomially, with independent $n-k$ trials, each having $n$ equally likely outcomes. The Poissonization
device yields that $\Bbb P(\bold X=\bold x) \le cn^{1/2}\Bbb P(\bold Z=\bold x)$, where $\bold Z=(Z_1,\dots,Z_n)$ and $Z_j$ are {\it independent\/} copies
of $\text{Poisson }(1-k/n)$. 

Introduce the probability measure $\Bbb P^*$ defined on the space $\mathcal S$ of triples $(\bold x, G, H)$ by
\begin{multline*}
\Bbb P^*\bigl(\{\bold X=x\}\cap\{\mathcal G=G\}\cap \{\mathcal H=H\}\bigr)\\
=\Bbb P(\bold Z=\bold x)\cdot \Bbb P\bigl(\{\mathcal G=G\}\cap\{\mathcal H=H\}
\boldsymbol|\,\bold X=\bold x\bigr).
\end{multline*}
Then $\Bbb P (E)\le cn^{1/2}\, \Bbb P^*(E)$ for all $E\subseteq \mathcal E$. By switching to $\Bbb P^*$ we gain independence of $X_1,\dots, X_n$
at the expense of the $cn^{1/2}$ factor. In particular, $\mathcal P_{nk}\le cn^{1/2}\Bbb P^*(A\cap B\cap C)$. So we turn to upper-bounding
{
$\Bbb P^*(A\cap B\cap C)$. To this end, we claim first that
\begin{equation}\label{2n}
\begin{aligned}
\Bbb P^*(B)&=(1-t)^{n-k+1}\left[1-f(t)p_m(t)\right]^{n-k+1} \\
p_m(t)&=\sum_{\ell=0}^m\frac{(1-t)^{\ell}}{\ell!},\,\, t:=\tfrac{k}{n},\,\,f(t):=t e^{-1+t}.
\end{aligned}
\end{equation}
}
 The first factor is the probability that none of columns from $[n-k+1]$ selects a row from $[k]$ in the first round, and the second factor is the
 probability no column $j\in [n-k+1]$, such that $X_j\le m$ would select a vertex in $[k]$ in the second round. Here we used the independence
 of $X_j$ under $\Bbb P^*$ and 
 \[
 \Bbb P^*(X_j\le m)=\sum_{\ell=0}^m e^{-1+t} \frac{(1-t)^{\ell}}{\ell!},\quad t=\frac{k}{n}.
 \]

So we need to estimate $P_{nk}=\Bbb P^*(A\cap C\,\boldsymbol|\,B)$ the probability of $A\cap C$, conditioned on the event $B$: every column in $[n-k+1]$
selects a row in $[n-k]$, and---if it is unpopular among those rows---would select such a row again. Let $\mathcal S$ stand for the full description of
selections by rows from $[k]$ in both rounds, and by columns from $[k-1]$ in first round, {\it compatible\/} with $A\cap B$. 

Let us specify, in four items,  a generic value $T$ of $\mathcal T$, $\mathcal T$ being a partial description of $\mathcal S$: {\bf (1)\/} let $V\subseteq [k-1]$ be the set of columns from $[k-1]$ whose first round choice rows are in $[k]$; {\bf (2)\/}  let $U\subset [k]$ be the set of the rows each
selected by at least $m+1$ columns, and {\bf (3)\/} so that $W:=[k]\setminus U$ is the set of unpopular rows in $[k]$. Denote $u=|U|,\,v=|V|,\,w=|W|$; evidently $u(m+1)\le v\le k-1$, $w=k-u$. For $i\in [k]$, let $a_i$ be the number of columns from $V$ which selected row $i$. 
{\bf (4)\/} To finish description of $T$, for $j\in [k-1]$, let $b_j$ be the number of rows in $[k]$ whose first round selection is the column $j$, and let $\be_j$ be the number of unpopular rows, those from $W$, whose second round selection is column $j$. On the event $A$, we have $\sum_{j\in [k-1]} b_j=k$ and $\sum_{j\in [k-1]}\be_j=k-u$.

Let $p_j(b_j,\be_j)=\Bbb P(F_j)$, $F_j$ the event that column $j$ has at least two contacts with rows in $[k]$. For $j\not\in V$, we have
\begin{equation}\label{3n}
\begin{aligned}
p_j(b_j,\be_j)\!&=\!\Bbb I\{b_j+\be_j\ge 2\}+ \Bbb I\{b_j=1, \be_j=0\} f(t) p_{m-1}(t)\\
&\quad +\Bbb I\{b_j=0, \be_j=1\}p_m(t)f(t),\quad f(t):=te^{-1+t},\,\,\,t=k/n.
\end{aligned}
\end{equation}
{\it Explanation.\/} For $j\notin V$, the first choice of column $j$ is a row in $[n-k]$. Suppose $b_j=1,\,\be_j=0$. $F_j$ holds if $X_j\le m-1$, making $j$
unpopular and allowing $j$ a second round selection, that happens to be a row in $[k]$, an event of probability 
\[
\left(e^{-1+t}\sum_{\ell=0}^{m-1} (1-t)^{\ell}/\ell!\right)\cdot t=p_{m-1}(t) f(t).
\]
 Suppose $b_j=0,\,\be_j=1$. 
This time $F_j$ holds if $X_j\le m$ and, again, $j$'s second selection is a row in $[k]$, an event of probability $p_m(t) f(t)$.

For $j\in V$, ($j$'s first choice is a row in $[k]$),  the counterpart of \eqref{3n} is
\begin{equation}\label{4n}
p_j(b_j,\be_j)=\Bbb I\{b_j+\be_j\ge 1\}+\Bbb I\{b_j+\be_j=0\}p_m(t)f(t).
\end{equation}

Conditioned on $\{\mathcal T(\mathcal S)=T\}\cap B$, the events $F_j$ are independent, so that
\[
\Bbb P^*(C\,\big|\,B\cap \{\mathcal T(\mathcal S)=T\}\big)=\! P^*\!\Biggl(\bigcap_{j-1}^kF_j\,\Big|\,B\cap\{\mathcal T(\mathcal S)=T\}\!\!\Biggr)=\prod_{j=1}^k p_j(b_j,\be_j).
\]
The RHS is the explicit function of $\mathcal T=\mathcal T(\mathcal S)$, and we need to compute its expected value to obtain $P_{nk}=\Bbb P^*(A\cap C|\,B)$. Denoting $\binom{\mu}{\vec{\boldsymbol\nu}}=\mu!/\vec{\boldsymbol\nu}!=\mu!/(\nu_1!\cdots\nu_t!)$, $(\nu_1+\cdots+\nu_t=\mu)$,
the resulting formula is
\begin{equation}\label{5n}
\begin{aligned}
&\qquad\qquad\,\, P_{nk}=\sum_{0\le u\le v\le k-1}\!\!P_{nk}(u,v),\\
&\qquad P_{nk}(u,v)=\binom{k-1}{v}\left(\frac{n-k}{n}\right)^{k-1-v}\!\!\!\sum_{\sum_i a_i=v\atop |\{i: a_i>m\}|=u}\!\!\binom{v}{\vec{\bold a}} n^{-v}\\ 
&\qquad\times \sum_{\sum_jb_j=k,\,\sum_j\be_j=k-u}
\!\binom{k}{\vec{\bold b}} n^{-k} \,\binom{k-u}{\vec{\boldsymbol\be}} n^{-k+u}\prod_j p_j(b_j,\be_j).
\end{aligned}
\end{equation}
(In the last product we can assume that $V=[v]$.) Given $u$ and $v$, the range of $\vec{\bold a}$ is non-empty only if $u(m+1)\le v$.

{\it Explanation.\/} We (1) select $v$ columns from $[k-1]$, and note that the probability that each of the remaining $k-1-v$ columns selects a row in
$[n-k]$ is $\left(\frac{n-k}{n}\right)^{k-1-v}$ ; (2) partition the chosen $v$ columns into an
ordered sequence of $k$ sets of cardinalities $a_1,\dots, a_k\ge 0$, with exactly $u$ $a_i$s above $m$, so that $a_i$ columns  select 
row $i\in [k]$, with overall probability $n^{-v}$;  (3) partition rows in $[k]$ (the set of $(k-u)$ unpopular rows in $[k]$, i.e. those chosen by  at most $m$
columns resp.) into an ordered sequence of subsets of cardinalities $b_1,\dots,b_{k-1}$ ($\be_1,\dots,\be_{k-1}$ resp.), so that each column $j$
is selected by $b_j$ rows in the first round  (by $\be_j$ unpopular rows in the second round resp.), with overall probability $n^{-k}\cdot n^{-(k-u)}$; (4)
add the contributions coming from each triple of partitions weighted with the factors $\prod_j p_j(b_j,\be_j)$.

To find a tractable upper bound for $P_{nk}(u,v)$ we will use the generating functions and the Chernoff-type bound: for non-negative
sequence $\{g_{\mu}\}$, and $x>0$, we have $g_m\le x^{-m} g(x)$, where $g(x)=\sum_{\mu\ge 0} g_{\mu}x^{\mu}$; analogous inequality holds for
multivariate generating functions with non-negative coefficients. Needless to say, this approach is contingent
on availability of an explicit formula for $g(x)$. 

The bottom sum does not depend on $\vec{\bold a}$. 
By symmetry, we  have  
\begin{equation}\label{6n}
\begin{aligned}
&\qquad\qquad\qquad\sum_{\sum_i a_i=v\atop |\{i: a_i>m\}|=u}\!\!\binom{v}{\vec{\bold a}} =v!\binom{k}{u}\sum_{a_1,\dots, a_u>m\atop
a_{u+1},\dots, a_k\le m}\frac{1}{\vec{\bold a}!}\\
&=v!\binom{k}{u}\,[x^v] \left(\sum_{a>m}\frac{x^a}{a!}\right)^u\cdot \left(\sum_{a\le m}\frac{x^a}{a!}\right)^{k-u}
\!\!\!\!\le v!\binom{k}{u} x^{-v} \exp^u_{m+1}(x) \cdot q^{k-u}_m(x),\\
&\qquad\qquad\quad\qquad{\exp_s(x):=\sum_{\tau\ge s}x^{\tau}/\tau!,\, q_s(x):=\sum_{\tau\le s} x^{\tau}/\tau!. }
\end{aligned}
\end{equation}
 Similarly, by \eqref{3n}-\eqref{4n},  the bottom sum in \eqref{5n} equals
\begin{multline*}
\sum_{\sum_jb_j=k,\,\sum_j\be_j=k-u}
\binom{k}{\vec{\bold b}}\,\binom{k-u}{\vec{\boldsymbol\be}} \prod_j p_j(b_j,\be_j)\\
=k!\,(k-u)!\, [y^k z^{k-u}]\prod_j\left(\sum_{b, \be}\frac{y^b z^{\be}}{b! \be !} p_j(b, \be)\right).
\end{multline*}
Denoting $\eta:=y+z$, the last product equals 
\begin{align*}
&\Biggl(\sum_{b,\be} \frac{y^b z^{\be}}{b! \be !} \Bigl(\!\Bbb I\{b_j+\be_j\ge 2\}\!+\! \Bbb I\{b_j=1, \!\be_j=0\} p_{m-1}(t)f(t)\\
&\qquad\qquad\qquad\qquad\qquad\qquad\qquad +\!\Bbb I\{b_j=0,\! \be_j=1\}p_m(t)f(t)\!\Bigr)\!\!\Biggr)^{k-1-v}\\
&\qquad\qquad\qquad\times\left(\sum_{b,\be} \frac{y^b z^{\be}}{b! \be !} \Bigl(\Bbb I\{b_j+\be_j\ge 1\}+\Bbb I\{b_j+\be_j=0\}p_m(t)f(t)\Bigr)\!\!\right)^{v}\\
&\qquad\qquad\qquad\qquad\quad=\left(\sum_{s\ge 2}\frac{(y+z)^s}{s!}+yp_{m-1}(t)f(t)+zp_m(t)f(t)\!\!\right)^{k-1-v}\\
&\qquad\qquad\qquad\qquad\qquad\qquad\quad\times \left(\sum_{s\ge 1}\frac{(y+z)^s}{s!}+p_m(t)f(t)\!\!\right)^{v}\\
&\quad\quad =\Bigl(\exp_1(\eta)+p_m(t)f(t)\Bigr)^{v}\cdot \Bigl(\exp_2(\eta)+yp_{m-1}(t)f(t)+zp_m(t)f(t)\Bigr)^{k-1-v}.
\end{align*}
Therefore we have
\begin{multline}\label{7n}
\sum_{\sum_jb_j=k,\,\sum_j\be_j=k-u}
\binom{k}{\vec{\bold b}}\,\binom{k-u}{\vec{\boldsymbol\be}} \prod_j p_j(b_j,\be_j)\\ \le \frac{k! (k-u)!}{y^k z^{k-u}} 
\Bigl(\exp_1(\eta)+p_m(t)f(t)\Bigr)^{v}\\
\times \Bigl(\exp_2(\eta)+yp_{m-1}(t)f(t)+z p_m(t)f(t)\Bigr)^{k-1-v}.
\end{multline}
Combining \eqref{5n}, \eqref{6n} and \eqref{7n}, {\it and\/} recalling the bound $u(m+1)\le v$, we conclude that
{
\begin{align}\label{8n}
P_{nk}&\le \frac{(k!)^2{(1-t)^{k-1}}}{n^{2k}y^k z^k}\,q^k_m(x)\Bigl(\exp_2(\eta)+(yp_{m-1}(t)+zp_m(t))f(t)\Bigr)^{k-1}\\
\notag
&\qquad\qquad \times \sum_{0\le u(m+1)\le v\le k-1}\!\!\!\!\frac{(k-1)_v }{u!} \,\xi^u \zeta^ v,\quad 
\xi:=\frac{nz\exp_{m+1}(x)}{q_m(x)},\\
\notag
\zeta&:={ \frac{g}{nx(1-t)}},\,\,\,g=g(t;y,z):=\frac{\exp_1(\eta)+p_m(t)f(t)}{\exp_2(\eta)+(yp_{m-1}(t)+zp_m(t))f(t)}.
\end{align}
}
Crucially, the sequence of the sums in \eqref{8n} has a simple (exponential) generating function. Indeed if $\zeta w<1$, then we have
\begin{align*}
&\sum_{k\ge 1}\frac{w^{k-1}}{(k-1)!}\left(\sum_{0\le u(m+1)\le v\le k-1}\frac{(k-1)_v}{u!}\,\xi^u\zeta^v\right)\\
&=\sum_{0\le u(m+1)\le v}\frac{\xi^u\zeta^v}{u!}\sum_{k-1\ge v}\frac{w^{k-1}(k-1)_v}{(k-1)!}
= e^w \sum_{0\le u(m+1) \le v}\!\!\frac{\xi^u\zeta^vw^v}{u!}\\
&=e^w\sum_{u\ge 0}\frac{(\xi(\zeta w)^{m+1})^u}{u!}\sum_{j\ge 0}(\zeta w)^j=\frac{\exp\bigl(w+\xi(\zeta w)^{m+1}\bigr)}{1-\zeta w}.
\end{align*}
Therefore we obtain: for all { $\bold R=(x,\,y,\,z,\,w) >\bold 0$}, such that $w<\zeta^{-1}$,
\begin{equation}\label{9n}
\begin{aligned}
&P_{nk}\le Q_{nk}(\bold R)\!:=\!\frac{(k!)^2{(1-t)^{k-1}}\Bigl(\exp_2(\eta)+(yp_{m-1}(t)+zp_m(t))f(t)\Bigr)^{k-1}}{n^{2k}y^k z^k}\\
&\qquad\times \,\frac{(k-1)! \,q^k_m(x)}{w^{k-1}}\cdot \frac{\exp\Bigl(w+\xi(\zeta w)^{m+1}\Bigr)}{1-\zeta w}; \qquad (\eta=y+z).
\end{aligned}
\end{equation}
At the price of $(x,y,z,w)$, yet to be chosen, we got rid of the multi-fold summation. Denoting $\rho=w/n$, we have
$\zeta w={g\rho/[(1-t)x]}$. The first line expression and the second line expression in \eqref{9n} are respectively of orders
\begin{equation*}
\begin{aligned}
&\frac{kg}{\exp_1(\eta)+p_m(t)f(t)}\cdot{(1-t)^k}\left(\frac{k}{ne}\right)^{2k}\left(\frac{\exp_1(\eta)+p_m(t)f(t)}{yzg}\right)^k;\\
&\quad\frac{k^{1/2}\rho}{t} \left(\frac{tq_m(x)}{e\rho}\right)^k\cdot\frac{\exp\left[n\left(\rho+\frac{z\exp_{m+1}(x)}{q_m(x)}\cdot\left({\frac{g\rho}{(1-t)x}}\right)^{m+1}\right)\right]}{1-{\frac{g\rho}{(1-t)x}}},\\
\end{aligned}
\end{equation*}
{\it uniformly\/} for all admissible { $\bold R$}; we used \eqref{8n} for both
expressions. 
So, denoting $\bold r=(x,y,z, \rho)$,  \eqref{9n} becomes
\begin{equation}\label{10n}
\begin{aligned}
P_{nk} =& O\!\left(\frac{k^{3/2}\rho g}{t\Bigl(\exp_1(\eta)+p_m(t)f(t)\Bigr)}\cdot \exp\bigl(n\mathcal H_{n,m}(t;\bold r)\bigr)\!\!\right),\\
\mathcal H_{n,m}(t;\bold r):=&2t\log\frac{t}{e} +{t\log(1-t)+}\,
 t\log\left(\frac{\exp_2(\eta)+(yp_{m-1}(t)+zp_m(t))f(t)}{yz}\right) \\
&\qquad + t\log\left(\frac{tq_m(x)}{e\rho}\right) +
  \rho+
  \frac{z\exp_{m+1}(x)}{q_m(x)}\cdot\left({\frac{g\rho}{x(1-t)}}\right)^{m+1} \\
&\qquad - n^{-1}\log\left(1-{\frac{g\rho}{x(1-t)}}\right).
\end{aligned}
\end{equation}
Recall that $P_{nk}=\Bbb P^*(A\cap C| B)$ and $\Bbb P^*(B)$ is given by \eqref{2n}. Therefore
\begin{equation}\label{11n}
\mathcal P_{nk}\!=\!O\bigl(n^{1/2}P_{nk}\,\Bbb P^*(B)\!\bigr)\!=O\!\left(\!n^{1/2}Q^*_{nk} \bigl[(1-t)(1-p_m(t) f(t))\bigr]^{n-k}\right),
\end{equation}
and, by \eqref{1n},
\begin{equation}\label{12n}
E_{nk}=O\left(\frac{k}{n}\binom{n}{k}^2\mathcal P_{nk}\!\right).
\end{equation}
Collecting the estimates \eqref{10n}-\eqref{12n},  
we arrive at
\begin{subequations}
\begin{align}
\label{13n}
E_{nk} =& O\left(\frac{(nk)^{1/2}\rho g}{t\bigl(\exp_1(\eta)+p_m(t)f(t)\bigr)}\cdot\exp\bigl(nH_{n,m}(t;\bold r\bigr)\!\right),\\
\label{13n2}
H_{n,m}(t;\bold r):=& -2t+(1-t)\log\frac{1-p_m(t)f(t)}{1-t}\, {+ t \log(1-t)}\\
\notag
&\hspace{-.5in} + t\log\left(\frac{\exp_2(\eta)+(yp_{m-1}(t)+zp_m(t))f(t)}{yz}\right) +
  t\log\left(\frac{tq_m(x)}{e\rho}\right)\\
\notag
&\hspace{-.5in} +\rho+\frac{z\exp_{m+1}(x)}{q_m(x)}\cdot\left({\frac{g\rho}{x(1-t)}}\right)^{m+1}-n^{-1}\log\left(1-{\frac{g\rho}{x(1-t)}}\right).
\end{align}
\end{subequations}
where the $H_{n,m}$ in equation~\eqref{13n2} is not the same as the 
$\mathcal H_{n,m}$ in~\eqref{10n} because of the inclusion of terms from
equations~\eqref{1n} and \eqref{2n}.   We already proved (see \eqref{2n++}) that $\sum_{k\le (m+1-\eps)\log n}E_{nk}=O(n^{-\eps/(m+1)})$, if $m>0$. So the task is to establish existence of the tuple $\mathbf r$ for every $t\ge (m+1-\eps)\frac{\log n}{n}$ such that $nH_{n,m}(t;\mathbf{r})$ is negative enough to out-power the front factor in \eqref{13n}, so that the sum of $E_{n,k}$ over the remaining $k$'s will go to zero as well. It is beneficial to start
earlier, with $t\ge 1/n$, i.e. with $k\ge 1$.

 \subsection{Small $t$'s.}
 Our focus is on $m>0$, but for comparison we include here $m=0$ as well. {Intuitively, for small $t$'s the search for the sub-optimal $\mathbf{r}$
ought to be done by narrowing down the field of candidates. After some tinkering with $H_{n,m}(t;\mathbf r)$ we chose  $\bold r(t)=(x(t),y(t), z(t),\rho(t))$: }
\begin{equation}\label{xyzrho}
x(t)=a t^{\sigma},\quad y(t)=b_1 t^{\sigma},\quad z(t)=b_2 t^{\sigma}, \quad \rho(t)=c\,t,
\end{equation}
with the parameters to be determined. Then, calculating upper bounds for the various terms in \eqref{13n2}, 
\begin{align*}
&-2t+(1-t)\log\frac{1-p_m(t)f(t)}{1-t}\,{ +t\log(1-t)}=-t(1+e^{-1}p_m(0)) +O(t^2);\\
&\begin{aligned}
  t\log\left(\frac{\exp_2(\eta)+\bigl(yp_{m-1}(t)+zp_m(t)\bigr)f(t)}{yz}\right)&=t\log\frac{\frac{(b_1+b_2)^2t}{2} +O(t^{1+\sigma})}{b_1b_2t}\\
&=t\log \frac{(b_1+b)^2}{2b_1b_2} +O(t^{1+\sigma});
\end{aligned} \\
&t\log\left(\frac{tq_m(x)}{e\rho}\right)+\rho=t\log \frac{1+O(t^{\sigma})}{ec} +ct= t\left(\log\frac{1}{ec}+c\right) +O(t^{1+\sigma});\\
&\frac{z\exp_{m+1}(x)}{q_m(x)}\cdot\left({\frac{g\rho}{x(1-t)}}\right)^{m+1} =\left\{\begin{aligned}
&\frac{2b_2c}{b_1+b_2}t +O(t^{1+\sigma}),&& m=0,\\
&O(t^{2-\sigma}),&& m\ge 1;\end{aligned}\right.\\
&1-\frac{g\rho}{x(1-t)}= 1-(1+O(t^{\sigma}))\frac{2ct^{1-2\sigma}}{a(b_1+b_2)}=1-\Theta\bigl(t^{1-2\sigma}\bigr);
\end{align*}
{ since $1-g\rho/[x(1-t)]$ is the argument of {the} $\log$-function in \eqref{13n2}, we need to choose $\sigma<1/2$.}
{Combining} the estimates we obtain: 
\begin{align*}
H_{n,m}(t;\bold r(t))&=-\gamma_m (\bold b,c) t+O(t^{1+\sigma}) +\Theta\bigl(n^{-1} t^{1-2\sigma}\bigr),\\
\gamma_m (\bold b,c)&=1+e^{-1}p_m(0)-\log\frac{(b_1+b_2)^2}{2b_1b_2}-\left(\log\frac{1}{ec}+c\right)\\
&\quad -\frac{2b_2c}{b_1+b_2}\cdot\Bbb I(m=0).
\end{align*}
It follows that $H_{n,m}(t;\bold r(t))$ is continuous, but {\it not differentiable\/} at $t=0$.
So $\gamma_m (\bold b,c)$ depends on $c$ and $b_2/b_1$ only, while $\sigma$ determines the behavior of the remainders.  With some calculus it follows that, for $m>0$, $\gamma_m (\bold b,c)$ attains its
maximum at $b_2/b_1=c=1$, while for $m=0$ the maximum is attained at $b_2/b_1=c=1/\sqrt{3}$. Explicitly,
\begin{equation}\label{14-n}
\gamma_m:=\max\gamma_m(\bold b,c)=\left\{\begin{aligned}
&1+e^{-1}-\log\bigl(\sqrt{3}+2\bigr),&& m=0,\\
&1+e^{-1}\sum_{j\le m} 1/j! - \log 2,&& m>0;\end{aligned}\right.
\end{equation}
$\gamma_0=0.0509\dots$, and $\gamma_1=1.0426\dots$, {\it exceeding $\gamma_0$ by a $20+$ factor\/}, while $\gamma_{\infty}=1.3068\dots$.
(In this regard the case $m=0$ is drastically different from  the case $m>0$.) Therefore
\begin{equation}\label{14n}
H_{n,m}(t)=\min_{\bold r} H_{n,m}(t;\bold r)\le -\gamma_m t +O(t^{1+\sigma})+\Theta\bigl(n^{-1}t^{1-2\sigma}\bigr).
\end{equation}
Consequently, for every $m\ge 0$, and $n$ sufficiently large, the function $H_{n,m}(t)$ is negative for $t\in [n^{-\la}, \eps_m]$, where $\la\in \bigl(1, (2\sigma)^{-1}\bigr)$, and $\eps_m>0$ is chosen sufficiently small.  
 
The front factor by $\exp\bigl(nH_{n,m}(t;\bold r\bigr)$ in \eqref{13n} is of order $(nk)^{1/2}$.
So it follows from \eqref{13n} and \eqref{14n} that for $m\ge 1$
 \begin{equation}\label{16n}
 E_{n.k} =O\left((nk)^{1/2}\exp\Bigl(-\gamma_m k+O\bigl(n^{-1/2}k^{3/2}\bigr)\Bigr)\right).
 \end{equation}
 We proved already that, for $k=O(n^{1/2})$ and $m\ge 1$,
\[
E_{nk}\le n^{-1}\exp\Bigl(\frac{k}{m+1}+O\bigl(k^{\frac{m}{m+1}}+k^2/n\bigr)\Bigr).
\]
Therefore
\begin{align}
\log E_{nk}&\le \min\Bigl\{\!-\log n+k/(m+1)),\,0.5\log n -\gamma_m k\Bigr\}\label{17.5n}\\
&\quad +O\Bigl(n^{-1/2}k^{3/2}+k^{\frac{m}{m+1}}+\log k\Bigr).\notag
\end{align}
The explicit term \eqref{17.5n}  attains its {\it maximum\/} at 
\[
k_n=\frac{1.5(m+1)}{(m+1)(1+e^{-1}p_m(0)-\log 2)+1}\cdot \log n
\]
and the maximum is $-c_m \log n$,
$
c_m:=1-\frac{1.5}{1+(m+1)\gamma_m}.
$
$c_m$ increases with $m$, $c_1\approx 0.514$ and $c_{\infty}=1$.  It follows easily that for $m\ge 1$
\begin{equation}\label{18n}
\max\{\log E_{nk}:\,k\le n^{1/2}\} \le -c_m\log n +O\Bigl((\log n)^{\frac{m}{m+1}}\Bigr),
\end{equation}
so that
\[
\sum_{k\in [1,k_n]}E_{nk}=O\Bigl(k_n\exp(-c_m\log n +O((\log n)^{\frac{m}{m+1}})\Bigr)=n^{-c_m +o(1)}.
\]
And
\begin{align*}
\sum_{k\in [k_n, n^{1/2}]}E_{nk}&=O\left(\exp\Bigl(-c_m\log n +O((\log n)^{\frac{m}{m+1}}\Bigr)\sum_{j\ge 0}e^{-(\gamma_m+o(1))j}\right)\\
&=n^{-c_m+o(1)}.
\end{align*}
We conclude that, for $m\ge 1$, $\sum_{k\le n^{1/2}}E_{nk}\le n^{-c_m+o(1)}$.

\subsection{Moderate and large $t$'s.} At this final stage we concentrate on $m=1$. Plugging $m=1$ into the formula for $H_{n,m}(t;\bold r)$ in
\eqref{13n2} we obtain
 \begin{align*}
&H_{n,1}(t;\bold r=(x,y,z,\rho)):=-2t+(1-t)\log\frac{1-(2-t)f(t)}{1-t}\,{+t\log(1-t)}\\
&\qquad+t\log\left(\frac{e^{\eta}-1-\eta+\bigl(y+z(2-t)\bigr)f(t)}{yz}\right)+t\log\left(\frac{t(1+x)}{e\rho}\right)\\
&\qquad+\rho+\frac{z(e^x-1-x)}{1+x}\cdot\left({\frac{g\rho}{x(1-t)}}\right)^2-n^{-1}\log\left(1-{\frac{g\rho}{x(1-t)}}\right);\\
&g=\frac{e^{\eta}-1+(2-t)f(t)}{e^{\eta}-1-\eta+\bigl(y+z(2-t)\bigr)f(t)},\quad f(t)=te^{-1+t},\quad \eta=y+z.
\end{align*}
In light of the analysis in the previous section, we need to check that
$H_{n,1}(t):=\min_{\,\bold r}H_{n,1}(t;\bold r)<0$ for $t\in [n^{-1/2},1/2]$ and $n$
sufficiently large.
Extensive numerical analysis involving {three} independent Matlab
minimization algorithms demonstrated that 
\[
\max \{ H_{n,1}(t) : t \in [1/n,1/2]\} = H_{n,1}(1/n) < 0 \text{ for all }n \geq 100.
\]
In particular, $H_{n,1}(1/2)\approx -0.051$ and $\min H_{n,1}(t)\approx -0.065$. That the function $H_{n,1}(t)$ is ``barely negative'' for $t\in (0,1/2]$
may be charitably interpreted as supporting our decision to use Matlab software, rather than to search for a protracted, yet uninspiring, calculus-based
proof of $H_{n,1}(t)$'s negativity. See Appendix for the details. 
\qed
\\

Though negative, $\gamma_0$ in \eqref{14-n} is so close to zero, that one wonders whether $H_{n,0}(t)$ will remain negative for all $t\in (0,1/2]$, like
all other $H_{n,m}(t)$ for $m>0$.  As discussed  in Appendix, once $t$ exceeds $0.0035$, $H_{n,0}(t)$ becomes and remains
positive for all remaining $t\leq1/2$; see the figure.  
This figure shows the striking difference between $H_{n,m}$ for $m = 0$ and~1.
\begin{figure}[ht]
\begin{center}
\includegraphics[scale=0.66,clip]{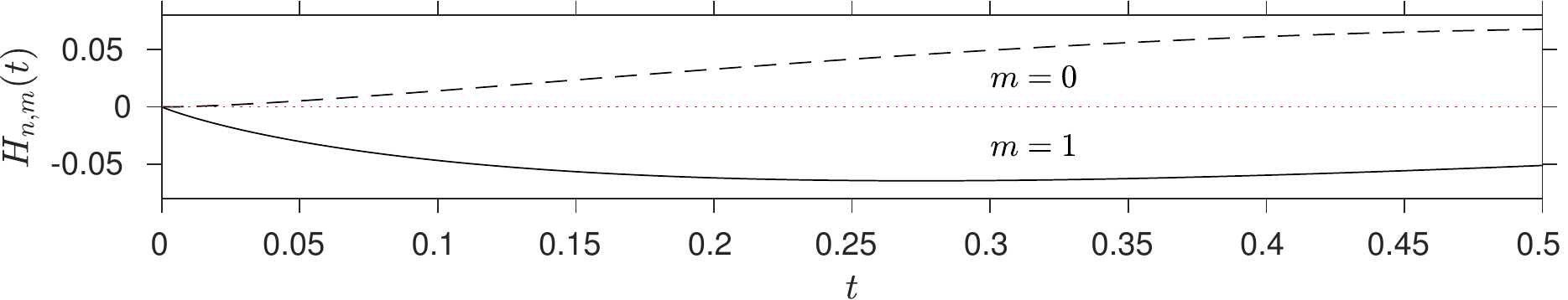}
\caption{\label{fig}
  The scaled logarithm of the expected number of obstacles to a perfect
  matching, $H_{n,m}(t)$ with $n = 100 \cdot 2^{10}$, 
  for $m = 0$ and 1.}
\end{center}
\end{figure}
At the very least, it means that, if true, almost sure existence of a perfect matching in $B_{n,0}$ would require a different argument. More plausibly
though, we need to consider a possibility that existence of a perfect matching in $B_{n,0}$ is unlikely. In the next Section we show that this is indeed the case. The proof itself will not require help of Matlab.\\ 

{\bf Part 2.\/} Let us prove that for $m>0$
\[
\Bbb P_{n,m}:=\Bbb P(B_{n,m}\text{ is connected})= 1-O\Bigl(n^{-c_m+o(1)}\Bigr).
\]
On the event $\mathcal E_n:= \{B_{n,m} \text{ is not connected and has a perfect matching}\}$, there exist a row set $K$ and a column set $L$ such that
$|K|=|L|=k$, $k\le n/2$, that induce a component of $B_{n,m}$. Let $X_k$ denote the total number of such pairs.
Denote $k_1=k$, $k_2=n-k$. We have 
\begin{align*}
\Bbb E[X_k]&\le \binom{n}{k}^{\!2}\prod_{i=1}^2\left(\frac{k_i}{n}\right)^{2\left(2k_i-\frac{k_i}{m+1}\right)}.
\end{align*}
{\it Explanation.\/} Let $K_1=K$, $L_1=L$, $K_2=V_1\setminus K_1$, $L_2=V_2\setminus L_1$. On the event ``$K,L$ induce a component of $B_{n,m}$''
in round $1$ every column from $L_i$ selects a row in $K_i$ and every row from $K_i$
selects a column in $L_i$. The total number of unpopular columns in  $L_i$ is { at least} $k_i-\frac{k_i}{m+1}$; so the probability that the selections
by columns from $L_i$ in the two rounds are all among rows in $K_i$ is  { at most} $(k_i/n)^{2k_i-\frac{k_i}{m+1}}$. Likewise 
$(k_i/n)^{2k_i-\frac{k_i}{m+1}}$ is { an upper bound for the}  probability that the selections by rows from $K_i$ in the two rounds are all among columns
in $L_i$.

Using $\binom{n}{k}\le n^n/[k^k (n-k)^{n-k}]$, we obtain then
\begin{align*}
\Bbb E[X_k]&\le \binom{n}{k}^{\! 2}\left(\frac{n^n}{k^k(n-k)^{n-k}}\right)^{-\frac{4m+2}{m+1}}\le \binom{n}{k}^{\!-\frac{2m}{m+1}}.
\end{align*}
Consequently
\[
\Bbb P(\mathcal E_n)\le \sum_{k=1}^{n/2}\Bbb E[X_k]=O\left[\binom{n}{1}^{-\frac{2m}{m+1}}\right]=O\bigl(n^{-\frac{2m}{m+1}}\bigr),
\]
implying that
\begin{align*}
\Bbb P(B_{n,m}\text{ is not connected}) &=O\bigl(n^{-\frac{2m}{m+1}}\bigr)+\Bbb P(B_{n,m}\text{ has no perfect matching})\\
&=O\bigl(n^{-\frac{2m}{m+1}}\bigr) +O(n^{-c_m+o(1)})=O(n^{-c_m+o(1)}).
\end{align*}
\qed

\section{Proof of Theorem \ref{2}.\/}

Let $k=\lfloor n^{\delta}\rfloor$, $\delta\in (0,1)$ to be specified shortly,  
and let $Y_n$ be the number of set pairs
$(K,L)$, $|K|=k$, $|L|=k-1$ such that $L=\Gamma(K)$. This time 
we {\it drop\/} the condition that every column vertex in $L$ has at least two row neighbors in $K$. It suffices to show that a.a.s.\@ $Y_n>0$;
indeed by Hall's Marriage Lemma existence of such a pair rules out existence of a perfect matching.
To this end we will prove that $\Bbb E[Y_n]\to\infty$ and $\Bbb E[Y_n^2]\sim \Bbb E^2[Y_n]$.  

{\bf (1)\/} For $0\le u\le v\le k-1$, The counterpart of $P_{nk}(u,v)$ in \eqref{5n} is given by 
\[
P_{nk}(u,v)= \frac{(k-1)^{2k-u}(k)_u}{n^{2k-u+v}}\cdot \binom{k-1}{v}\left(1-\frac{k}{n}\right)^{k-1-v}\!\!S(v,u);
\]
here $S(v,u)$ is the Stirling number of the second kind, i.e. the number of partitions of $[v]$ into $u$ non-empty sets.                                     

\noindent {\it Explanation.\/} We (1) choose $v$ columns from $[k-1]$ and $u$ rows from $[k]$ in $\binom{k-1}{v}\binom{k}{u}$ ways; (2)
allocate $v$ columns among $u$ rows in $u! S(v,u)$ ways, so in round $1$ the remaining $k-1-v$ columns select rows from $[n-k]$; (3) allocate $k$ rows among $k-1$ columns, thus determining round $1$ selections of columns
in $L$ made by rows from $[k]$, and allocate $k-u$ unpopular rows among $k-1$ columns, thus determining round $2$ selections of columns
still in $L$ made by unpopular rows from $[k]$, in $(k-1)^{k+(k-u)}$ ways total. Finally
\begin{equation}\label{19n-}
\frac{1}{n^{2k-u+v}}\left(1-\frac{k}{n}\right)^{k-1-v}
\end{equation}
is the probability of each of the resulting outcomes. We need a sharp asymptotic formula for $P_{nk}:=\sum_{0\le u\le v\le k-1} P_{nk}(u,v)$. Notice
that for $u<v$, by log-concavity of $\{S(v,u)\}_{u\le v}$ (Harper \cite{Har}, Godsil \cite{God}, Section 6.3) we have
\begin{align*}
\frac{P_{nk}(u+1,v)}{P_{nk}(u,v)}&=\frac{n(k-u)}{k-1}\cdot\frac{S(v,u+1)}{S(v,u)}\\
&\ge \frac{n(k-u)}{k-1}\cdot\frac{S(v,v)}{S(v,v-1)}\ge \frac{n}{k}\cdot \frac{1}{v}\ge \frac{n}{k^2}\ge \frac{n}{n^{2\delta}}\to\infty,
\end{align*}
if $\delta<1/2$. So we have
\[
\sum_{u=0}^v P_{nk}(u,v)=\bigl(1+O(n^{-1+2\delta})\bigr) P_{nk}(v,v),
\]
 with the front factor absorbing $\left(1-\frac{k}{n}\right)^{k-1-v}$ from \eqref{19n-}, which implies
\begin{equation}\label{19n}
\begin{aligned}
P_{nk}&=\sum_{0\le u\le v\le k-1}P_{nk}(u,v)=\bigl(1+O(n^{-1+2\delta})\bigr)\frac{(k-1)^{2k}}{n^{2k}}S_k,\\
&\qquad\quad S_k:=\sum_{v\le k-1}\frac{(k)_v}{(k-1)^v}\binom{k-1}{v}.
\end{aligned}
\end{equation}
The ratio of two consecutive terms in the sum $S_k$ decreases with $v$. So the largest term corresponds to the smallest $v$ for which this ratio is 
at most one. It easily follows this $v$ is one of two integers closest to $v_0=\sigma k+(-2 +2/\sqrt{5}),\quad \sigma:=\frac{3-\sqrt{5}}{2}$.
The dominant contribution to $S_k$ comes from the terms with $|v-v_0|\le k^{1/2}\log k$, and uniformly for these $v$ we have
\begin{align*}
&\frac{(k)_v}{(k-1)^v}\binom{k-1}{v}=(1+O(k^{-1}))\frac{e^{\sigma}}{\sqrt{2\pi\sigma k}}\,e^{kH(v/k)},\\
&\quad H(z):=-z-2(1-z)\log (1-z)-z\log z.
\end{align*}
Not surprisingly, $H(z)$ attains its maximum at $z=\sigma$. Approximating $H(v/k)$ by the quadratic Taylor polynomial around $\sigma$, we replace
the sum by the Gaussian integral and obtain
\begin{align*}
&\quad S_k=(1+O(k^{-1}))\frac{\exp(\sigma+kH(\sigma))}{\sqrt{\sigma(-H''(\sigma))}},\\
&H(\sigma)=-\sigma-2\log(1-\sigma)=0.5804576362.
\end{align*}
This formula for $S_k$ results in a compact estimate of $P_{nk}$ in \eqref{19n}. We hasten to add that we have not considered yet another condition
a pair $(K,L)$ needs to meet:  (1) in round $1$ every column from $[n-k+1]$ selects a row in $[n-k]$;  (2) in round $2$ every column in $[n-k+1]$, which was not selected by any row from $[n-k]$ in round $1$, still selects a row from $[n-k]$. The event (1) has probability $(1-k/n)^{n-k+1}=\exp(-k+O(k^2/n))$; the event (2) is dependent on the event whose probability $P_{nk}$ we have analyzed. Its conditional probability 
is $(1-k/n)^{W}$, where $W$ is the number of columns in $[n-k+1]$, which are unpopular among rows in $[n-k]$. Since $k\ll n$, i.e. $k=o(n)$, by the Poissonization approximation we have: for $\eps\in (0,1/2)$,
\begin{equation}\label{20n-}
\Bbb P\bigl\{|W-e^{-1}n|\le n^{1/2+\eps}\bigr\}\ge 1-\exp\bigl(-\Theta( n^{2\eps})\bigr).
\end{equation}
And we observe that the total number of all pairs $(K,L)$,  with $|K|=k=\lfloor n^{\delta}\rfloor$,  is 
\[
\binom{n}{k}\binom{n}{k-1}=e^{O(k\log n)}=e^{O(n^{\delta}\log n)}\ll e^{\Theta(n^{2\eps})},
\]
provided that $\delta<2\eps$. Therefore, computing the moments of $Y_n$, we can---at the cost of an additive
error term $e^{-\Theta( n^{2\eps})}$---replace $W$ with $e^{-1}n+O\bigl(n^{1/2+\eps}\bigr)$, in which case $(1-k/n)^W= \exp\bigl(-e^{-1}k +O(n^{-1/2+\delta+\eps})\bigr)$. In combination with \eqref{19n} it follows then: for $0<\delta<\min(1/2-\eps,\,2\eps)$,
\begin{equation}\label{20n}
\begin{aligned}
\Bbb E[Y_n]&=\bigl(1+O(n^{-1/2+\delta+\eps})\bigr)\binom{n}{k}\binom{n}{k-1} e^{-k(1+e^{-1})}P_{nk} +O\bigl(e^{-\Theta( n^{2\eps})}\bigr)\\
&=\bigl(1+O(n^{-1/2+\delta+\eps})\bigr)\binom{n}{k}\binom{n}{k-1}e^{-k(1+e^{-1})}\frac{(k-1)^{2k}}{n^{2k}}S_k \\
&\null\qquad +O\bigl(e^{-\Theta( n^{2\eps})}\bigr)\\
&=\bigl(1+O(n^{-1/2+\delta+\eps})\bigr)\frac{(k-1)^{2k}}{k! (k-1)!\,n} e^{-k(1+e^{-1})}S_k+O\bigl(e^{-\Theta( n^{2\eps})}\bigr)\\
&=\Theta\bigl(n^{-1}e^{\la k}\bigr),\quad\la:=1-e^{-1}+H(\sigma)=1.212578195 >0.
\end{aligned}
\end{equation}

{\bf (2)\/}  Next we will show that $\Bbb E[(Y_n)_2]\lesssim \Bbb E^2[Y_n]$, where $(Y_n)_2=Y_n(Y_n-1)$ is the total number of ways to select two bad pairs,
 $(K_1,L_1)$ and $(K_2,L_2)$, $|K_i|=k$, $|L_i|=k-1$. Given $0\le \mu\le k$, $0\le \nu\le k-1$, every $\{(K_i,L_i)\}_{i=1,2}$ with $|K_1\cap K_2|=\mu$, 
 $|L_1\cap L_2|=\nu$ has the same probability, call it $\Pi(\mu,\nu)$, that $(K_1,L_1)$ and $(K_2,L_2)$ are both bad. The contribution of all such $(\mu,\nu)$ configurations to $\Bbb E[(Y_n)_2]$ is 
 \[
 \Pi(\mu,\nu)\binom{n}{k-\mu,\,\mu,\,k-\mu}\binom{n}{k-1-\nu,\,\nu,\,k-1-\nu}.
 \]
 So we focus on $K_1=[k]$, $L_1=[k-1]$, $K_2=[k-\mu+1, 2k-\mu]$, $L_2=[k-\nu, 2k-2-\nu]$. Visually, we have two $k\times (k-1)$ rectangles on the $n\times n$ integer lattice,
 the first rectangle occupying the North-West corner, and the second rectangle having its North-West corner point at $(k-\mu+1, k-\nu)$.
 
Let us define $P_{nk}(\bold u,\bold v)$, a  counterpart of $P_{nk}(u,v)$. { Here $\bold v=(v_1, v_2, \hat v_1, \hat v_2)$, $\bold u=(u_1, u_2)$,
\[
v_i\le |L_i\setminus(L_1\cap L_2)|,\,\, \hat v_1+\hat v_2\le |L_1\cap L_2|,\,\,u_i\le |K_i\setminus(K_1\cap K_2)|.
\]
Introduce the event $A(\bold u,\bold v)$:  $v_i$ ($\hat v_i$ resp.) is the number of columns belonging to $L_i\setminus(L_1\cap L_2)$ ($L_1\cap L_2$ resp.) that in round $1$ selected a row belonging  only to $K_i$, and $u_i$ is the number of those rows. $P_{nk}(\bold u,\bold v)$ is the probability of the event $A(\bold u,\bold v)$ } {\it intersected\/} with the event that in
round $1$ no row from $K_i$ selected a column from $L_i^c$, and in round $2$ no unpopular row belonging only to $K_i$ selected a column from 
$L_i^c$. $P_{nk}(\bold u,\bold v)$ is an {\it upper\/}  bound for the probability that both $(K_1,L_1)$ and $(K_2, L_2)$ are bad, and the two probabilities
are equal when $\mu=\nu=0$. Arguing like $P_{nk}(u,v)$, we have
\begin{multline*}
P_{nk}(\bold u,\bold v)=\binom{\nu}{\hat v_1,\hat v_2}\left(1-\frac{2(k-\mu)}{n}\right)^{\nu-\hat v_1-\hat v_2}\\
\times\prod_{i=1}^2\binom{k-1-\nu}{v_i}\binom{k-\mu}{u_i} u_i!\, S(v_i+\hat v_i,u_i)(k-1)^{2k-\mu-u_i}\\
\times\left(1-\frac{k-\mu}{n}\right)^{k-1-\nu-v_i}\frac{1}{n^{2k-u_i+v_i+\hat v_i}}.
\end{multline*}
 {\it Explanation.\/} The trinomial coefficient and the four binomial coefficients should be clear. $u_i!\, S(v_i+\hat v_i, u_i)$ is the number of ways to
 assign $v_i+\hat v_i$ columns from $L_i$  to the already chosen $u_i$ row vertices belonging only to $K_i$, with each of these rows getting at least
 one column. $(k-1)^{2k-\mu-u_i}$ is the total number of ways to assign $k$ rows from $K_i$ and $k-u_i$ unpopular
 rows belonging only to $K_i$ to the columns from $L_i$.  $\bigl(1-2(k-\mu)/n\bigr)^{\nu-\hat v_1-\hat v_2}$ is the probability that, by the definition of $\hat v_1$ and $\hat v_2$, some specific $\nu-\hat v_1-\hat v_2$ common columns each have to select a row outside of the symmetric difference $K_1\Delta K_2$.
$\bigl(1-(k-\mu)/n\bigr)^{k-1-\nu-v_i}$ is the probability that the columns from $L_i\setminus (L_1\cap L_2)$ select a row outside 
$K_i\setminus (K_1\cap K_2)$. $1/n^{2k-u_i+v_i+\hat v_i}$ is the probability that both first round
selections by rows from $K_i$ and second round selections by unpopular rows belonging to $K_i\setminus (K_1\cap K_2)$ are among columns from $L_i$. 

To estimate $\sum_{\bold u,\bold v}P_{nk}(\bold u, \bold v)$  for guidance we use elements of the part {\bf (1)\/}.  First of all, both
$\bigl(1-2(k-\mu)/n\bigr)^{\nu-\hat v_1-\hat v_2}$ and $\bigl(1-(k-\mu)/n\bigr)^{k-1-\nu-v_i}$ are each equal to $1+O(n^{-1+2\delta})$, since
$k=\Theta(n^{\delta})$.  Second,
by log-concavity of {$(k-\mu)_u\,S(v_i+\hat v_i, u)$} (as a function of $u$) we have
\begin{multline*}
\sum_{u_i\le v_i+\hat v_i} \binom{k-\mu}{u_i} u_i! S(v_i+\hat v_i,u_i)\frac{(k-1)^{2k-\mu-u_i}}{n^{2k-u_i+v_i+\hat v_i}}\\
\le\bigl(1+O(n^{-1+2\delta})\bigr) (k-\mu)_{v_i+\hat v_i}\frac{(k-1)^{2k-\mu-v_i-\hat v_i}}{n^{2k}};
\end{multline*}
{the front factor on the RHS is that close to $1$ because the ratio of the last term and the penultimate term of the sum is at least $n/(k-1)^2$.}
So
\begin{equation}\label{22n}
\begin{aligned}
\sum_{\bold u,\bold v}P_{nk}(\bold u, \bold v)&{ \le }\bigl(1+O(n^{-1+2\delta})\bigr)\\
&\qquad\times \sum_{v_1,v_2}\,\,\prod_{i=1}^2\binom{k-1-\nu}{v_i}(k-\mu)_{v_i}
\frac{(k-1)^{2k-\mu-v_i}}{n^{2k}}\\
&\qquad\qquad\times\sum_{\hat v_1,\hat v_2}\binom{\nu}{\hat v_1,\hat v_2}\prod_{i=1}^2 (k-\mu-v_i)_{\hat v_i} (k-1)^{-\hat v_i};
\end{aligned}
\end{equation}
$v_i\le k-1-\nu,\quad \hat v_1+\hat v_2=:s\le \nu$. The bottom sum in \eqref{22n} is
\begin{equation*}
\begin{aligned}
&\sum_{s\le \nu}\frac{\nu! (k-1)^{-s}}{(\nu-s)!}\!\!\sum_{\hat v_1+\hat v_2=s}\prod_{i=1}^2\!\binom{\!k-\mu-v_i}{\hat v_i\!}\\
&=\sum_{s\le \nu}\frac{\nu! (k-1)^{-s}}{(\nu-s)!\!}\binom{2(k-\mu)-v}{s}\le \sum_{s\le \nu}\binom{\nu}{s} (k-1)^{-s}(2k)^s\\
&=\left(\frac{2k}{k-1}+1\right)^{\nu}\le 4^{\nu}. 
\end{aligned}
\end{equation*}
The penultimate sum in \eqref{22n} is at most
\[
 \left(\frac{(k-1)^{2k}}{n^{2k}}\right)^2\sum_{v_1,v_2}\prod_{i=1}^2\frac{(k)_{v_i}}{(k-1)^{v_i}}\binom{k-1}{v_i}.
\]
Therefore the equation \eqref{22n} becomes
 \begin{multline}\label{23n}
\sum_{\bold u,\bold v}P_{nk}(\bold u, \bold v)\le\bigl(1+O(n^{-1+2\delta})\bigr) 4^{\nu} (k-1)^{-2\mu}\\
\times \left(\frac{(k-1)^{2k}}{n^{2k}}\right)^2\Biggl[\sum_w\frac{(k)_w}{(k-1)^w}\,\binom{k-1}{w}\Biggr]^2\\
=\bigl(1+O(n^{-1+2\delta})\bigr)\,4^{\nu}\, (k-1)^{-2\mu}\Biggl(\frac{(k-1)^{2k}}{n^{2k}}\,S_k\Biggr)^2,
\end{multline}
{according to the definition of $S_k$ in \eqref{19n}.}
The LHS sum is the probability that in round $1$ rows  from $K_i$ select columns from $L_i$, and that in round $2$ the unpopular rows belonging
to $K_i\setminus (K_1\cap K_2)$ again each select column in $L_i$.  $\Pi(\mu,\nu)$ is, {\it at most\/}, the probability of this event intersected with the event:  everyone of $n-(2(k-1)-\nu)$ columns in $(L_1\cup L_2)^c$ selects in round $1$ one of $n-(2k-\mu)$ rows in $(K_1\cup K_2)^c$, and each of the $W_1$ columns unpopular among these rows
selects a row in $(K_1\cup K_2)^c$ in round $2$ again.  Since $k=\Theta(n^{\delta})$, analogously to $W$ we have
\[
\Bbb P\bigl\{|W_1-e^{-1}n|\le n^{1/2+\eps}\bigr\}\ge 1-\exp\bigl(-\Theta(n^{2\eps})\bigr),
\]
see  \eqref{20n-}. Therefore with probability that high, the conditional probability of the event ``$W_1$ columns stay with rows from $(K_1\cup K_2)^c$
in round $2$'' is at most 
\[
\left(\!1-\frac{2k-\mu}{n}\right)^{n-(2(k-1)-\nu)} \!\left(1-\frac{k}{n}\right)^{W_1}\!\!=\bigl(1+O(n^{-1+2\delta})\bigr) e^{-2k(1+e^{-1})+e^{-1}\mu}.
\]
Therefore, by \eqref{23n}
\begin{align*}
\Pi(\mu,\nu)&\le \bigl(1+O(n^{-1+2\delta})\bigr)\,4^{\nu}\, (k-1)^{-2\mu}\Biggl(\frac{(k-1)^{2k}}{n^{2k}}\,S_k\Biggr)^2\\
&\qquad\times \exp\bigl(-2k(1+e^{-1})+e^{-1}\mu\bigr).
\end{align*}
So we have
\begin{multline*}
\Bbb E[Y_n(Y_n-1)]\le O\bigl(e^{-\Theta(n^{2\eps})}\bigr)+ \bigl(1+O(n^{-1+2\delta})\bigr) e^{-2k(1+e^{-1})}\Biggl(\frac{(k-1)^{2k}}{n^{2k}}\,S_k\Biggr)^2\\
\times \sum_{\mu,\nu\ge 0} \binom{n}{k-\mu,\mu,k-\mu}\binom{n}{k-1-\nu,\nu,k-1-\nu} 4^{\nu} \left(\frac{e^{e^{-1}}}{(k-1)^2}\right)^{\mu}.
\end{multline*}
Both $\left\{\binom{n}{k-\mu,\,\mu,\, k-\mu}\right\}$ and $\left\{\binom{n}{k-1-\nu,\,\nu,\,k-1-\nu}\right\}$ are log-concave as functions of $\mu$ and $\nu$,
respectively. So the sum is at most 
\begin{align*}
&\binom{n}{k,0,k}\binom{n}{k-1,0,k-1}\sum_{\mu\le k}\left(\frac{k^2}{n-2k+1}\,\frac{e^{e^{-1}}}{(k-1)^2}\!\right)^{\mu}\cdot \sum_{\nu\le k-1}\!\!\left(\frac{4(k-1)^2}{n-2k+3}\right)^{\nu}\\
&\qquad = \bigl(1+O(k^2/n)\bigr) \frac{n^{2k}}{(k!)^2}\frac{n^{2(k-1)}}{((k-1)!)^2}=\bigl(1+O(k^2/n)\bigr) \binom{n}{k}^2\binom{n}{k-1}^2.
\end{align*}
Combining this equation and \eqref{20n} (second line), and recalling that $k=\Theta(n^{\delta})$, we obtain
\[
\Bbb E[Y_n(Y_n-1)]=O\bigl(e^{-\Theta(n^{2\eps})}\bigr)+\bigl(1+O(n^{-1/2+\delta+\eps})\bigr) \Bbb E^2[Y_n],
\] 
implying that for $0<\delta<\min(1/2-\eps, 2\eps)$,
\[
\frac{\Bbb E[Y_n^2]}{\Bbb E^2[Y_n]}=1+O(n^{-1/2+\delta+\eps}),
\]
since $\Bbb E[Y_n]=\exp\bigl(\Theta(n^{\delta})\bigr)$, see \eqref{20n}. By Chebyshev's inequality, 
\[
\Bbb P\bigl (Y_n\ge 0.5 \Bbb E[Y_n]\bigr)\ge 1- O(n^{-1/2+\delta+\eps}).
\]
\qed

\section{Components of $B_{n,0}$.} In \cite{KarPit} it was asserted that a.a.s.\@ $B_{n,0}$ consists of a single giant component and 
small isolated cycles (cyclic components) with a bounded total size. The proof was based on observation that in presence of a perfect matching
every isolated component $(K,L)$ must be balanced, i.e. $|K|=|L|$. However we know now that a.a.s.\@ $B_{n,0}$ has no perfect matching. Here is
a sketch of the corrected proof of a close claim; the only computer aid it relies on is a surface plot.

Suppose that a pair $(K,L)$, ($|K|=k$, $|L|=\ell$), induces a component of $B_{n,0}$. We focus on {\it smaller\/} components, i.e. of size $k+\ell\le n$. Introduce $(K_1,L_1)=(K,L)$, $(K_2,L_2)= (V_1\setminus K, V_2\setminus L)$, $k_1=k$, $k_2=n-k$, $\ell_1=\ell$, $\ell_2=n-\ell$. Suppose $\ell\le k$, i.e.
$\ell_1\le k_1$; then $\ell_2\ge k_2$. Let us bound the probability $P_{k,\ell}$ that
none of the pairs $(i,j)$ with $i\in K_1,\,j\in L_2$ or with $i\in K_2,\,j\in L_1$ is an edge of $B_{n,0}$. We have
\begin{align*}
P_{k,\ell}&\le P_{k,\ell}^*:=\left(\frac{k_1}{n}\right)^{\ell_1}
\left(\frac{\ell_2}{n}\right)^{k_2}\cdot\,\,\left(\frac{\ell_1}{n}\right)^{2k_1-\ell_1}\,\cdot\left(\frac{k_2}{n}\right)^{2\ell_2-k_2}\\
&\quad\times c_1\sqrt{k_1}\left(1-e^{-\frac{k_1}{\ell_1}}\frac{k_2}{n}\right)^{\ell_1}\cdot \,c_2\sqrt{\ell_2}\left(1-e^{-\frac{\ell_2}{k_2}}\frac{\ell_1}{n}\right)^{k_2}.
\end{align*}
{\it Explanation.\/} $1$-st line: first factor is the probability that in round $1$ vertices in $L_1$ and $K_2$ select, exclusively from their {\it larger\/} partner sets $K_1$ and $L_2$; $2$-nd factor ($3$-rd factor resp.) is an upper bound for the probability that all vertices in $K_1$
and all unpopular vertices in $K_1$ (all vertices in $L_2$
and all unpopular vertices in $L_2$ resp.) select vertices from $L_1$ ($K_2$ resp.). $2$-nd line:  $1$-st factor is an upper bound for the probability that none of the unpopular vertices in $L_1$ selects a vertex from $K_2$ in round $2$; $2$-nd factor is an upper bound for the probability that
none of the unpopular vertices in $K_2$ selects a vertex from $L_1$  in round $2$. For instance, the first bound comes from approximating the
numbers of vertices in $K_1$, which selected 
the vertices from $L_1$ in round $1$, by the $\ell_1$-long sequence of independent Poissons,
each with parameter $k_1/\ell_1$.

So, denoting the expected number of such pairs $(K,L)$ by $E_{k,\ell}$, we have 
\begin{align*}
E_{k,\ell}&\le \binom{n}{k}\binom{n}{\ell}P_{k,\ell}^*\le c n^{1/2} \exp\bigl(n \Bbb H(k/n,\ell/n)\bigr),\\
\Bbb H(x,y)&:\!=-x\log x-(1-x)\log(1-x)-y\log y-(1-y)\log(1-y)\\
&\quad\,\,+y\log x+(1-x)\log(1-y)\\
&\quad\,\, +(2x-y)\log y+(1+x-2y)\log(1-x)\\
&\quad\,\, +y\log\left(1-e^{-\frac{x}{y}}(1-x)\right)+(1-x)\log\left(1-e^{-\frac{1-y}{1-x}}y\right).
\end{align*}
Since $k+\ell\le n$, $\ell\le k$, we are interested at $y\le x$, $x+y\le
1$. The 3D plot of $\Bbb H(x,y)$ reveals that $\Bbb H(x,y)<0$ for all $x+y>0$ and
$\Bbb H(0+,0+)=0$, the latter seen directly from the formula for $\Bbb H(x,y)$. Setting $y=zx$, $z\in [0,1]$, we obtain: for $x$ small, 
\begin{align*}
 \Bbb H(x,y)&=(1-z) x\log x+x\Bigl(2(1-z)\log z+\!z\bigl(\log(1-e^{-1/z})\!-\!e^{-1}\bigr)\!\Bigr)\\
&\quad+O(x^2)\\
&\le x \sup_{z\in [0,1]}\Bigl(2(1-z)\log z+z\bigl(\log(1-e^{-1/z})-e^{-1}\bigr)\Bigr)+O(x^2)\\
&\le -0.648\, x+O(x^2).
\end{align*}
It follows that for $\ell\le k$, $\a>0$ and small $\eps>0$ 
\begin{multline*}
\sum_{\a\log n\le k+\ell\le n} E_{k,\ell}\le\sum_{0.5\a\log n\le k \le\eps n} E_{k,\ell}+\sum_{k\ge \eps n,\, k+\ell\le n} E_{k,\ell}\\
\le cn^{1/2}\sum_{k\ge 0.5\a\log n}\!\! k\exp\bigl(-k(0.648-O(\eps))\bigr)\\
 +O\left(n^{2.5}\exp\Bigl(n\max\{H(x,y): y\le x,\, x+y\in [\eps, 1]\}\Bigr)\right)\\
=O\left(n^{1/2}(\log n) n^{-0.5\a (0.648-O(\eps))}\right)\to 0,
\end{multline*}
if $\a> 1.55$ and $\eps>0$ is sufficiently small. Thus a.a.s.\@ all components smaller than the largest component must be of size $1.55 \log n$ at most.
The expected total size of such components is $\sum_{k+\ell\le 1.55 \log n} E_{k,\ell}$, which is easily seen to be of order $O(n^{1/2+o(1)})$. 

We conclude
that a.a.s.\@ $B_{n,0}$ consists of a single giant component and some components each of size $1.55\log n$ at most, whose total size is a.a.s.\@ of order $O(n^{1/2+o(1)})$.\\

{\bf Acknowledgment.\/} We are genuinely grateful to Michael Anastos and Alan
Frieze for closely reading the old paper coauthored by the two of us, and for
pinpointing a consequential oversight in one of the bounds. The full validation of the approach
adopted in that  paper required a substantial extension and diversification of the initial techniques.

\section{Appendix}

We explain how the numerical calculations were carried out in Matlab to
minimize $H_{n,m}(t;{\mathbf r})$ in \eqref{13n2}. To begin, we rewrite the
equations somewhat to explicitly show the independent variables in
$\mathbf r$.
Since $\eta = y + z$, we replace its one occurrence.
Next, $p_m(t) = q_m(1 - t)$, so we replace it.
Finally, we have to be careful of the last term in $H_{n,m}$ to make sure its
complicated argument is never non-positive, so we replace it by $u$, i.e.,
$-n^{-1} \log u$, and solve for $\rho$.
The independent variables are now ${\mathbf r} = (x,y,z,u)$.
However, we do not remove $\rho$ completely from the equation for $H$ since in
two of its four occurrences it is simpler to leave it in rather than replacing
it by a complicated function of $u$.
Also, from numerical evidence, $\rho$ is a much simpler function of $t$, so we
can more easily estimate the asymptotic behavior of the independent variables.

Combining these modifications, we obtain
\begin{subequations}
\begin{align}
\notag
f(t) &= t e^{t-1} \\
\notag
q_k(w) &= \sum_{j=0}^k \frac{w^j}{j!} \qquad \bigl(\text{in Matlab } \texttt{q(w,k)}\bigr) \\
\notag
\exp_k(w) &= e^w - q_{k-1}(w) \qquad \bigl(\text{in Matlab } \texttt{expq(w,k)}\bigr) \\
\notag
g_m(t; y, z) &= \frac{ \exp_1(y+z)+q_m(1-t)f(t) }
          { \exp_2(y+z) + \bigl[ y q_{m-1}(1-t) + z q_m(1-t) \bigr] f(t) } \\
\label{rho_m}
\rho_m(t; x, y, z, u) &= \frac{(1-t)(1-u)x}{g_m(t;y,z)}\\ 
\label{H_n,m}
H_{n,m}(t;x,y,z,u) &= -2t + (1 - t)\log\left(\frac{1 - q_m(1-t)f(t)}{1 - t} \right)\,+t\log(1-t) \\
\notag
  &\hspace{-0.6in} {} + t \log \left( \frac{ \exp_2(y+z) + \bigl[ y q_{m-1}(1-t) + z q_m(1-t) \bigr]f(t) }
      { y z } \right) \\
\notag
  &\hspace{-0.6in} {} + t \log \left( \frac{ t q_m(x) }{ e \rho_m(t;x,y,z,u) } \right) 
               + \rho_m(t;x,y,z,u)  \\
\notag
  &\hspace{-0.6in} {} + \frac{z \exp_{m+1}(x)}{q_m(x)} \, (1 - u)^{m+1} -
      \frac{1}{n} \log u .
\end{align}
\end{subequations}
where we have explicitly included all the arguments in each function.
In each numerical run we fix $m$ and $n$, we define the anonymous
functions exactly as written above, using precisely these arguments,
and, recalling that $t = k/n$, calculate the minimum of $H_{n,m}(t;\mathbf r)$ for
each $t_k = k/n$ where $k \in \bigl[1, \lceil n/2\rceil \bigr]$.
Additionally, there are constraints on the independent variables that
\begin{equation}
\label{constraints}
x,y,z \geq 0 \text{ and }u \in (0,1] \text{ for all } t \in [0,1/2].
\end{equation}
For each $t$, we denote the location of the minimum by $\overline{\mathbf r}$
and the minimum value itself by
$\overline H_{n,m}(t) = H_{n,m}(t;\overline{\mathbf r})$.

We used three independent iterative minimization functions,
\texttt{fminsearch}, \texttt{fminunc}, and \texttt{fmincon}, in Matlab; the
latter two are in the optimization toolbox (which costs extra). 
They delivered strikingly close trajectories for all $m \geq 0$ and $n$. 
In particular, for $m = 1$ and $n \geq 100$ 
the trajectories are strikingly close and negative for $t \geq 1/n$ and
$n \geq 100$.
This provides strong numerical evidence that the analytical minimum 
$\min_{\,\mathbf{r}}H_{n,m}(t;\mathbf r)$ is negative for $t\ge 1/n$.

We began with the first one which uses the Nelder-Mead simplex algorithm,
that does not require the function to be differentiable, but also does not
guarantee it converges to a minimum.
To obtain as much accuracy as possible and to try to prevent ``approximate''
minima, the function and optimality tolerances were set to $10^{-8}$.
However, it is an unconstrained minimization method.  
So, as it is commonly done, we added a penalty function, namely
\begin{equation}
\label{penaltyfn}
P \bigl(h(-x)x + h(-y)y + h(-z)z + h(-u)u \bigr)^2
\end{equation}
with $P = 10^4$, to~\eqref{H_n,m}.
$h$ is the Heaviside step function which ``nudges'' the iterates to stay in the
constraint region~\eqref{constraints} whenever any of the variables become negative. 
For each $t_k$, $k > 1$, the initial iterate is the solution at $t_{k-1}$. 
The reason we start at $t_1=1/n$, rather than at $t_0=0$, is that our
admittedly limited analysis of the asymptotic behavior of $H_{n,m}(t)$ as $t\downarrow 0$, see
\eqref{14n}, suggests strongly that the function is not differentiable at $t=0$.
Extensive numerical evidence suggests that the initial iterates
can be chosen at $t = t_1$ from~\eqref{xyzrho}:
we let $\sigma = 1/3$, and $a=b_2=1$, $b_1 = \sqrt{3}$, $c = 1/\sqrt{3}$ for
$m = 0$, while $a = b_1 =b_2= c=1$ for $m > 0$, where $u$ is obtained from
$\rho$ by using~\eqref{rho_m}. 

We are now ready to discuss the results, and we continue to focus on
\texttt{fminsearch}, discussing the differences with the other minimization
functions as we go along.
The curves $\overline{H}_{n,m}(t)$ for $m = 0$ and $m = 1$ with
$n = 100 \cdot 2^{10}$ are shown in Figure \ref{fig}.

First, we get the case $m = 0$ out of the way.
For $n \lessapprox 22000$, $\overline{H}_{n,0}(t) > 0$ for all $t > 0$.
However, for larger values of $n$, $\overline{H}_{n,0}(t) < 0$ for small $t$.
The values of $t$ at which $\overline{H}_{n,0}(t)$ becomes positive are
$t = 0.00215$ for $n = 10^5$, 
$0.003162$ for $n = 10^6$,
$0.0033659$ for $n = 10^8$, 
$0.003369374$ for $n = 10^9$, and
$0.0033698094$ for $n = 10^{10}$,
so the switch point on $t$-axis certainly seems to be approaching
a rather small value as $n\to\infty$.

Next, from numerical evidence for $m = 1$, the trajectory is negative
for all $t>0$, if  $n \geq 100$.
To see that the curves are converging, we show $\overline{H}_{n,1}(t)$ at $t = 0.01$
for $n = 100 \cdot 2^j$ where $j \in [0, 14]$:
\begin{align*}
& 0.0045632, \, -0.0055698, \, -0.0063425, \, -0.006925, \, -0.007355, \\
&\quad -0.0076648, \, -0.0078813, \, -0.0080275, \, -0.0081226, \, -0.0081822, \\
  &\quad -0.0082184, \, -0.0082397, \, -0.0082521, \, -0.0082591, \text{ and } -0.0082630.
\end{align*}
And we show it at $t = 0.5$:
\begin{align*}
& -0.0125880, \, -0.028543 , \, -0.038172, \, -0.043832, \,  -0.04709, \\
&\quad -0.048934, \, -0.049964, \, -0.050533, \, -0.050844, \, -0.051014, \\
  &\quad -0.051105, \, -0.051154, \,  -0.051181, \, -0.051195, \text{ and } -0.051202.
\end{align*}
Again, the numbers certainly seem to be decreasing to a limiting value
$< 0$.

The second minimization function we used is \texttt{fminunc}, which is also
unconstrained.
It is based on a quasi-Newton method, specifically the
Broyden-Fletcher-Goldfarb-Shanno algorithm with a cubic line search procedure,
where the gradient is approximated numerically, while the Hessian is
approximated by a secant-like method in higher dimensions.
Over the entire numerically calculated interval $t \in (0,1/2]$,
the curves generated by \texttt{fminsearch} and \texttt{fminunc} are negative
and differ by $<2\times 10^{-7}$.

The third minimization function is \texttt{fmincon}, which uses interior-point
optimization.
It is the only function which allows constraints, so no penalty function is
applied. 
However, the resulting curve rapidly oscillated for 
$t \lessapprox 10^{-3}$, repeatedly assuming positive values.
These oscillations continued for $t \lessapprox 5 {\times} 10^{-3}$ although
the curve remained negative, although for larger values of $t$ the difference
from \texttt{fminsearch}'s curve did fall below by $2 {\times} 10^{-7}$.
This curve cannot be accepted; so what could have gone wrong?

These large amplitude oscillations looked like a manifestation of a numerical
instability, which requires a technical explanation.
Minimization algorithms often have difficulties, much
more than zero-finding algorithms.
The latter only require the first derivative of the function, called the
Jacobian, calculated either analytically or numerically; the former
require the gradient, first derivatives, and also some approximation to the
Hessian, second derivatives, which introduces more errors.
Also, zero-finding is inherently more accurate because, even only considering
one dimension, finding the point where a curve passes through the $x$ axis is
much more accurate than finding where it attains a minimum.
(As a simple example, if $y = f(x)$ passes through the $x$ axis with slope
$s \neq 0$, a change in $y$ by $\delta y$ results in a change in $x$ by
$\delta x = \delta y/s$, while if $y$ has a minimum which behaves like
$a (x - \xi)^2$, a change in $y$ near the minimum by $\delta y$ results in a
change in $x$ by $\delta x = \sqrt{\delta y/a}$, i.e., $\delta y$ has an
exponent of $1/2$ rather than 1, so a small error in $\delta y$ results in a
much larger error in $\delta x$.)
It seems that, somehow, because of the numerical approximation to the
gradient, followed by a secant-like approximation to the Hessian, and in a
region where the valley surrounding the minimum was very shallow, a small
error in the solution at $t_j$, when used as the initial guess for $t_{j+1}$,
caused a larger error. 
This generated a feedback loop which finally died out at
$t \approx 5 {\times} 10^{-3}$.

To improve the accuracy of the calculations, we used alternate algorithms in
\texttt{fminunc}, a trust region algorithm,  and \texttt{fmincon},
a trust-region-reflective algorithm, both of which require the gradient of the  
function to be calculated analytically (not shown).
When these more accurate algorithms were used, these two curves were always
negative, and the differences between all three, i.e., including
\texttt{fminsearch}'s, were always $< 10^{-7}$.
We stated earlier that \texttt{fminsearch} was the most accurate of all the
algorithms. 
This claim is supported by numerically approximating the second derivatives of
all five curves using second-order centered differences.
By eye, the second derivative decreased monotonically from 104 to 0.18 over
the entire interval using \texttt{fminsearch}.
For the other two functions, without the analytical gradient, there were
fluctuations over much of the interval of magnitudes about 1000, while, with
the analytical gradient, there were only fluctuations for $t \ll 1$ with
magnitudes of 300 to 600. 

As another, rather strong, test of the accuracy of the code, the program was
only run for small $t$'s so that the slope of
$\overline{H}_{n,m}(t)$ at $t=t_1$ could be compared to \eqref{14-n}.
 A straight line was fit to the first 100 points using least squares.
The results for $m = 0$ and $n = 10^5$, $10^6$, $10^7$, $10^8$, $10^9$, and
$10^{10}$ are $-0.008904$, $-0.03394$, $-0.04374$, $-0.04942$, $-0.04950$, and
$-0.05027$ as compared to $-\gamma_0 = -0.0509$.
The same calculation for $m = 1$ produces
$-0.9501$, $-1.007$, $-1.028$, $-1.037$, $-1.040$, and $-1.041$ as compared to
$-\gamma_1 = -1.0426$.

In conclusion we note, for readers without access to Matlab, that no
modifications were required in the code to use Octave (a free software package
which is mostly compatible with Matlab) with \texttt{fminsearch}.

\textit{The Matlab code is (hopefully) accessible on the journal's website.}

We include a pseudocode showing the ``guts'' of the program.
Most of the code is taken up in calculating the various functions 
and generating the plots.

\def\CARET{\,\widehat{\phantom{o}}\,}
\def\nquad{\null\quad}
\noindent
\nquad $f \leftarrow (t) \;\cdots$ ; \\
\nquad $q \leftarrow (w,k) \;\cdots$ ; \\
\nquad $expq \leftarrow (w,k) \;\cdots$ ; \\
\nquad $g \leftarrow (t,y,z) \;\cdots$ ; \\
\nquad $\rho \leftarrow (t,x,y,z,u) \;\cdots$ ; \\
\nquad $H \leftarrow (t,x,y,z,u) \;\cdots \\
\nquad \quad\quad\quad\quad + P*(h(-x)*x + h(-y)*y + h(-z)*z + h(-u)*u)\CARET2$;\\
\nquad $dH \leftarrow (t,x,y,z,u) \;\cdots$ ; \qquad\qquad // array containing gradient of $H$\\
\nquad $it \leftarrow 0$; \\
\nquad \textbf{for} $t = 1/n$ to $1/2$ by $1/n$ \\
\nquad \quad $it \leftarrow it + 1$; \\
\nquad \quad \textbf{if} $t == 1/n$ \\
\nquad \quad\quad \textbf{if} $m == 0$ \\
\nquad \quad\quad\quad $x\_ic \leftarrow t\CARET(1/3); \; y\_ic \leftarrow \sqrt{3}\, t\CARET(1/3); \;
                       z\_ic \leftarrow t\CARET(1/3); \\
\nquad \quad\quad\quad \rho\_ic \leftarrow t/\sqrt{3}\,;$\\
\nquad \quad\quad \textbf{else} \\
\nquad \quad\quad\quad $x\_ic \leftarrow t\CARET(1/3); \; y\_ic \leftarrow t\CARET(1/3); \;
                z\_ic \leftarrow t\CARET(1/3); \; \rho\_ic \leftarrow t;$\\
\nquad \quad\quad \textbf{end} \\
\nquad \quad\quad $u\_ic \leftarrow 1 - g(t,y\_ic,z\_ic)*\rho\_ic/((1 - t)*x\_ic);$ \\
\nquad \quad \textbf{else} \\
\nquad \quad\quad $x\_ic \leftarrow x\_st(it-1); y\_ic \leftarrow y\_st(it-1); z\_ic \leftarrow z\_st(it-1);$ \\
\nquad \quad\quad $u\_ic \leftarrow u\_st(it-1);$ \\
\nquad \quad \textbf{end} \\
\nquad \quad $\{x\_st(it), y\_st(it), z\_st(it), u\_st(it)\} \\
\nquad \qquad\qquad \leftarrow  \text{minimization function}(H, dH, x\_ic, y\_ic, z\_ic, u\_ic)$;  \\
\nquad \textbf{end}

\end{document}